\input amstex
\documentstyle{amsppt}
\input amsppt.sty
\voffset=2truecm \hoffset=2truecm
 \pageno=1 \NoRunningHeads
\NoBlackBoxes

\def\p{\frak p}
\def\R{\widehat R}

\def\m{\frak m}
\def\q{\frak q}
\def\dim{\operatorname{dim}}
\def\co-supp{\operatorname{co-supp}}

\def\Ann{\operatorname{Ann}}
\def\Ass{\operatorname{Ass}}

\def\Usupp{\operatorname{Usupp}}
\def\lr{\longrightarrow}
\def\Att{\operatorname{Att}}
\def\Supp{\operatorname{Supp}}

\def\rad{\operatorname{rad}}

\def\docao{\operatorname{ht}}
\def\Spec{\operatorname{Spec}}
\def\N-dim{\operatorname{N-dim}}
\raggedbottom \normalbaselines \NoBlackBoxes
 \vskip 1cm
\topmatter
\title    Top local cohomology and the catenaricity of the unmixed
support of a finitely generated module
\endtitle
\endtopmatter
 \document
 \vskip 0.5cm
\centerline {NGUYEN TU CUONG$^a$, NGUYEN THI DUNG}
\smallskip
\centerline{ Institute of Mathematics}

\centerline{18 Hoang Quoc Viet Road, 10307 Hanoi, Vietnam}

\centerline{\it E-mail adress: ntcuong\@math.ac.vn}

\bigskip

\centerline{LE THANH NHAN$^b$}

\centerline{ Department of Mathematics}

\centerline{Thai Nguyen Pedagogical University, Thai Nguyen, Vietnam}

\centerline{\it E-mail adress: trtrnhan\@yahoo.com}

\vskip 1cm

\noindent { \bf Abstract.}{ \footnote{ 
 \hfill\break {\bf {AMS Classification }} 13D45, 13E10.
 \hfill\break $^a$Senior Associate Member of ICTP,
Trieste, Italy. \hfill\break $^b$Junior Associate Member of ICTP,
Trieste, Italy.}
 Let $(R,\m)$ be a Noetherian local
ring and $M$ a finitely generated $R-$module with $\dim M=d.$ This
paper is concerned with the following property for the top local
cohomology module $H^d_\m(M)$:
$$\Ann (0:_{H^d_\m(M)}\p)=\p\ \text{for all prime ideals}\
\p\supseteq\Ann H^d_\m(M).$$
In this paper we will show that this property is
equivalent to the catenaricity of the unmixed support $\Usupp M$ of $M$ which is defined by $\Usupp M=\Supp 
M/U_M(0)$, where $U_M(0)$ is the largest
submodule of $M$ of dimension less than $d.$   Some
characterizations of this property in terms of system of
parameters as well as the relation between the unmixed supports of $M$ and of the $\m$-adic completion $\widehat M$
 are given.
  \bigskip

 \noindent {\bf 1.  Introduction}
 \medskip

Throughout this paper, let $(R,\m)$ be a Noetherian local ring,
$M$ a finitely generated $R-$module with $\dim M=d,$ and $A$ an Artinian
$R-$module. For each ideal $I$ of $R$, we denote by $V(I)$ the set
of all prime ideals containing $I.$
\medskip

 An elementary property of finitely generated modules
  is that $\Ann  (M/\p M)=\p$ for all $\p\in
 V(\Ann  M).$ The dual question for Artinian modules is to ask whether
 $$\Ann (0:_A \p)=\p\ \text{for all}\ \p\in V(\Ann A).\tag *$$
 In case $R$ is complete with respect to the $\m-$adic topology,
the property (*) is satisfied  for all Artinian $R-$modules $A$
because of the Matlis duality between the
category of Noetherian $R-$modules and the category of Artinian
$R-$modules. Unfortunately the property (*) is not satisfied in
general.  For example, let $R$ be the Noetherian local domain of
dimension $2$ constructed by D. Ferrand and M. Raynaund [FR]  (see also [Na, App., Exam. 2]) such
that its $\m-$adic completion $\R$  has an associated prime
$\widehat\q$ of dimension $1.$ Then the Artinian $R-$module
$A=H^1_\m(R)$ does not satisfy the property (*), cf [CN1].
However, it seems to us that the property (*) is an important
property of Artinian modules. For example,  the property (*) is
closely related to some questions on dimension for Artinian
modules.  In [CN1], it is shown that $\N-dim A=\dim R/\Ann A$
provided $A$ satisfies the property (*), where $\N-dim A$ is the
Noetherian dimension of $A$ defined by Roberts [R] (see also
[K2]). Note that this equality does not hold in general. Concretely, with the Artinian $R-$module $A=H^1_\m(R)$ as above, $\N-dim
A=1<2=\dim R/\Ann A$ although this ring $R$ is catenary
and the top local cohomology module $H^2_\m(R)$ satisfies the property (*).
\medskip

The purpose of this paper is to study the property (*) for the top
 local cohomology Artinian module $H^d_\m(M)$, $\dim M=d,$ and
its applications.  We will show that, although $\N-dim H^d_\m(M)$
and $\dim R/\Ann H^d_\m(M)$ are always equal to $d,$ the property (*) is not  necessarily satisfied for $H^d_\m(M)$. Then we
find conditions such that $H^d_\m(M)$ satisfies the property (*). It seems surprising to us, that this property is
equivalent to some important properties of $M$. In particular,  the
property (*) is satisfied for $H^d_\m(M)$ if and only if  the
 unmixed support $\Usupp M=\Supp M/U_M(0)$ is
catenary, where $U_M(0)$ is the largest submodule of $M$ of
dimension less than $d.$ 

From now on,  we denote by $\widehat R$ (resp. $\widehat M$) the $\m-$adic completion of $R$ (resp. $M$) and $\Usupp_{\R} \widehat M$ the unmixed support of $\widehat M$ as an $\R-$module.  The main result of this paper is the following theorem.

\proclaim{Theorem }  The following statements are equivalent:
\medskip
(i) $\Ann (0:_{H^d_\m(M)}\p )=\p$ for all $\p\in V(\Ann
H^d_\m(M)).$
\smallskip
(ii) $\Usupp M$ is catenary.

\smallskip
 (iii) $\Usupp M= \{\widehat\p \cap R :\ \widehat \p\in\Usupp_{\R} \widehat M\}.$
\smallskip
(iv) For every sequence $x_1,\ldots ,x_d$ of elements in $\m$,
$(x_1,\ldots ,x_d)$ is a system of parameters of $H^d_\m(M)$ if
and only if it is a system of parameters of $M/U_M(0).$
\endproclaim
Here, the notion of system of parameters for Artinian modules is defined according to  Section 2.
\smallskip

As an immediate consequence of the above main theorem , we have the following characterization for the catenaricity of a Noetherian local domain.

\proclaim{Corollary} Suppose that   $(R,\m )$ is a  Noetherian local domain of dimension $d.$ Then  $R$ is catenary if and only if $H^d_\m(R)$  satisfies the property (*).
\endproclaim

\smallskip
This paper is divided into $4$ sections. In Section 2 we introduce
the property (*) for Artinian modules and recall some basic  facts that we need in the sequel. In the  last two
sections we present  the proof of the above main theorem.  The characterizations of the
property (*) for $H^d_\m(M)$ in terms of system of parameters and
the relation between two sets $\Usupp M$ and $\Usupp_{\R}
\widehat M$ (the part
(i)$\Leftrightarrow$(iii)$\Leftrightarrow$(iv) of the main theorem) are
proved in Section 3.   In Section 4, we prove the equivalence
between the property (*) for $H^d_\m(M)$ and the catenaricity of $\Usupp M$ (the part (i)$\Leftrightarrow$(ii) of the main theorem). We also examine some non-catenary domains to clarify the
results.
\bigskip

 \noindent {\bf 2. When is $\Ann (0:_A\p )=\p$ for all prime ideals $\p\supseteq\Ann A$?}

\medskip

 For each Noetherian $R-$module $M,$ it is clear that
 $\Ann (M/\p M)=\p$ for each prime ideal $\p\in V(\Ann  M).$
Therefore, for each Artinian $R-$module $A$, we consider the
following property:
$$\Ann (0:_A \p)=\p\ \text{for all prime ideals}\ \p\in
V(\Ann A).\tag *$$
\medskip
  As we mentioned in the introduction,  the property
 (*) is not satisfied for all Artinian modules.
 In this section, we give some conditions such that this property is satisfied.
  \medskip

  First we have the following result which  is proved in [CN1].
\smallskip

\proclaim{Proposition 2.1}   $A$ satisfies the property (*)  if one of the following conditions  holds.

(i) $R$ is complete with respect to the $\m -$adic topology.

(ii) $A$ contains a submodule which is isomorphic  to the
injective hull of $R/\m $.
 \endproclaim

 It should be mentioned that $\Supp M=\{\widehat\p \cap R\ :\ \widehat\p \in\Supp_{\R}\widehat
M\}$ for each finitely generated $R-$module $M.$  This means that $V(\Ann  M)=\{\widehat\p \cap R\ :\
\widehat\p \in V(\Ann_{\R}(\widehat M)\}.$ We also note that, for
each  Artinian $R-$module $A$, $A$ has a natural structure as an
Artinian $\R-$module (cf. [Sh]), and with this structure, a subset
of $A$ is an $R-$submodule of $A$ if and only if it is an
$\R-$submodule of $A.$ Therefore  it is natural to ask whether
$$V(\Ann A)=\{\widehat\p\cap R\ :\ \widehat\p\in
V(\Ann_{\R}A\}.$$ Below we show  that this equality holds if and only if $A$ satisfies  the property (*).  Recall that $A$ has a minimal
secondary representation
 $A=A_1+\ldots +A_n,$ where $A_i$ is $\p_i-$secondary, and the set
 $\{\p_1,\ldots ,\p_n\}$ does not depend on the choice of the
 minimal secondary representation of $A.$ This set is denoted by
 $\Att A$ and called {\it the set of attached prime ideals} of
 $A,$ cf. [Mac].  It follows by [Sh] that $$\Att A=\{\widehat\p\cap R\ :\ \widehat\p\in\Att_{\R}A\}.$$

\proclaim{Proposition 2.2}  The following conditions are
equivalent:

(i)  $A$ satisfies the property (*).

(ii) $V(\Ann  A)= \{\widehat\p \cap R :\ \widehat \p\in
V(\Ann_{\R}A)\}.$
\endproclaim

\demo{Proof} (i)$\Rightarrow$(ii).  Let $\widehat\p\in V(\Ann_{\R}A).$ Then $\widehat\p\supseteq\widehat\q$ for some minimal prime ideal $\widehat\q$ containing $\Ann_{\R}A.$ Note that any minimal prime ideal containing $\Ann_{\R}A$ is a minimal element of  $\Att_{\R} A,$ cf. [Mac]. Therefore  $\widehat\q\in\Att_{\R}A.$  So, $\widehat\q \cap R\in\Att A.$ Hence $\widehat\q\cap R\in V(\Ann A)$ and hence $\widehat\p\cap R\in V(\Ann A).$  Conversely, let $\p\in V(\Ann A).$ Then $\Ann (0:_A\p
)=\p$ by the hypothesis (i). Since $\p$ is a minimal prime ideal
containing $\Ann (0:_A\p ),$ it follows  that
$\p\in\Att (0:_A\p ).$ Therefore there exists a prime ideal
$\widehat\p\in\Att_{\R}(0:_A\p )$ such that $\widehat\p\cap R=\p
.$ Since $\widehat\p\in\Att_{\R}(0:_A\p )$, we have
$\widehat\p\supseteq \Ann_{\R}(0:_A\p ),$ and hence
$\widehat\p\in V(\Ann_{\R}A)$ with $\widehat\p \cap R=\p .$
\smallskip

\noindent (ii)$\Rightarrow$(i).  Let $\p\in V(\Ann A).$ By the
hypothesis (ii), there exists $\widehat\p\in V(\Ann_{\R}A)$ such
that $\widehat\p\cap R=\p.$ It follows by Matlis duality that
$\Ann_{\R}(0:_A\widehat\p )=\widehat\p .$ Therefore
$$\p\subseteq\Ann (0:_A\p )\subseteq
\Ann_{\R}(0:_A\widehat\p )\cap R=\widehat\p\cap R=\p .$$  Thus, $\Ann (0:_A \p)=\p$.\ \hfill\ \qed
\enddemo

 Roberts [R] introduced the concept of Krull dimension for Artinian modules.
  Kirby [K2] changed the terminology of Roberts and used the terminology of
 Noetherian dimension to avoid confusion
 with Krull dimension defined for finitely generated modules.
 In this paper we use the terminology of Kirby [K2].
 The {\it Noetherian dimension} of $A,$
 denoted by $\N-dim_RA,$ is defined inductively as follows:  when $A=0,$  put $\N-dim_R
A=-1.$ Then by induction, for an integer $d\geq 0,$ we put
$\N-dim_R A=d$ if $\N-dim_R A<d$ is false and for every ascending
sequence $A_0\subseteq A_1\subseteq \ldots $ of submodules of $A,$
there exists $n_0$ such that $\N-dim_R  (A_{n+1}/A_n)<d$ for all
$n>n_0$.
\smallskip

The following result gives some good  properties of Noetherian
dimension for Artinian modules which are in some sense dual to
that of Krull dimension for Noetherian modules, cf [R], [K1].

\proclaim{Lemma 2.3} $\ell(0:_A\m^n)$ is a polynomial for $n \gg
0,$ and
 $$\N-dim_R A=\deg \ell(0:_A\m^n)=
 \inf\{ t: \exists x_1,\ldots ,x_t\in\m : \ell (0:_A(x_1,\ldots ,x_t)R)<\infty \}.$$
\endproclaim
It follows by Lemma 2.3 that there exists a sequence $(x_1,\ldots
,x_d)$, $d=\N-dim_R A,$ such that $\ell (0:_A(x_1,\ldots
,x_d)R)<\infty .$ A such sequence is called a {\it system of
parameters of $A$}.
\smallskip

Note that, with the natural structure as an Artinian $\R-$module,
$\N-dim_RA=\N-dim_{\R}A.$ Therefore there is no confusion in
writing $\N-dim A$ instead of $\N-dim_RA$ or $\N-dim_{\R}A.$
Moreover, it follows by Matlis duality that
$$\N-dim A=\dim\widehat R/\Ann_{\widehat R}A
=\max\{\dim \R/\widehat\p \ :\ \widehat \p\in\Att_{\R}A\}.$$
 Note that the set of all minimal elements of $\Att A $ is exactly the set of all minimal prime
ideals containing $\Ann A$. So, we have
$$\dim R/\Ann A=\max\{\dim R/\p\ :\ \p\in\Att A \}.$$
The following result shows  the relation between  $\N-dim A$ and $\dim
R/\Ann A.$  

\proclaim{Proposition 2.4} [CN1]. The following statements are
true.

(i) $\N-dim A\leqslant \dim R/\Ann A.$

(ii) If  $A$ satisfies the property (*) then 
$\N-dim A=\dim R/\Ann A.$
\endproclaim
\medskip

\noindent {\bf Remark 2.5.}

 (i) As we mentioned in the introduction, there exist Artinian modules $A$ such that $\N-dim A<\dim R/\Ann A$.
\smallskip

(ii) The converse of Proposition 2.4,(ii) is not true. In the next
sections, we will show that the top local cohomology module
$A=H^d_\m(M)$ of a finitely generated $R-$module $M$  of dimension $d$ does not
satisfy the property (*) in general, although it always satisfies
the condition $\N-dim A=\dim R/\Ann A=d.$

\bigskip

\noindent {\bf 3. The property (*) for the top local cohomology
modules}

 \medskip

From now on, let  $M$ be a finitely generated $R-$module with
$\dim M=d.$ In this section, we examine the property (*)  for the
 top local cohomology  module $H^d_\m (M).$
\medskip
We first present a nice property of $H^d_\m (M)$, cf. [CN1,
Corollary 3.6].

  \proclaim{Lemma 3.1}  $\N-dim H^d_\m(M)=\dim R/\Ann H^d_\m(M)=d.$
  \endproclaim

Let $U_M(0)$ be the largest submodule of $M$ of dimension less
than $d.$  Note that if $\displaystyle
0=\bigcap_{\p\in \Ass M}N(\p )$ is a reduced primary decomposition
of the zero submodule of $M$ then $U_M(0)=\displaystyle \bigcap_{\dim
R/\p =d} N(\p ),$ cf. [CN2].  Therefore we have
$$\Ass M/U_M(0)=\{\p\in \Ass M\ :\ \dim R/\p =d\}.$$ Hence
$$\Supp M/U_M(0)=\bigcup_{\p\in\Ass M,\ \dim
R/\p =d}V(\p ).$$ The set $\Supp M/U_M(0)$ is called {\it the
unmixed } support of $M$ and denoted by $\Usupp M.$ 

\proclaim{Lemma 3.2} Let $\p\in\Supp M.$ Then $\p\in\Usupp M$ if
and only if $\p\supseteq\Ann H^d_\m(M).$ In particular, $\Usupp
M=V(\Ann H^d_\m(M)).$
\endproclaim
\demo{Proof} We have by [BS] that 
$$\Att H^d_\m(M)=\{\q\in\Ass
M\ :\ \dim R/\q =d\}.$$ 
Moreover,  the set of all minimal prime ideals containing $\Ann H^d_\m(M)$ and the
set of all minimal elements of $\Att H^d_\m(M)$ are the same.
Therefore
$$V(\Ann H^d_\m(M))=\bigcup_{\p\in\Ass M,\ \dim
R/\p =d}V(\p )=\Usupp M.$$ \ \hfill\ \qed
\enddemo

There are some nice relations between associated primes and the supports of $M$ and of its $\m-$adic completion $\widehat M.$  For example,  $\Ass M=\{\widehat\p\cap R\ :\ \widehat\p\in\Ass_{\R}\widehat M\}$ and  $\Supp M=\{\widehat\p\cap R\ :\ \widehat\p\in\Supp_{\R}\widehat M\}.$ Moreover,
$$\{\p\in\Ass M\ :\ \dim R/\p =d\}=\{\widehat\p\cap R\ : \widehat\p\in\Ass_{\R}M, \  \dim \R/\widehat\p =d\}.$$
So, it is natural to ask about the relation between $\Usupp M$ and $\Usupp_{\R}\widehat M.$  First we have the following lemma.

\proclaim{Lemma 3.3}  $\Usupp M\supseteq \{\widehat\p \cap R :\ \widehat
\p\in\Usupp_{\R} \widehat M\}.$
\endproclaim
\demo{Proof} Let $ \widehat
\p\in\Usupp \widehat M.$ Then $ \widehat \p\supseteq \widehat\q$ for some $\widehat \q\in \Ass_{\R} \widehat M$ satisfying $\dim \R/\widehat\q =d.$  It follows that $\widehat \q\cap R\in \Ass M$ and $\dim R/(\widehat \q\cap R)=d.$  Since $ \widehat \p\cap R\supseteq  \widehat \q\cap R,$ we get  $ \widehat \p\cap R\in \Usupp M.$ \ \hfill\ \qed
\enddemo
 In general,  the two sets  $\Usupp M$ and $ \{\widehat\p \cap R :\ \widehat
\p\in\Usupp_{\R} \widehat M\}$ are different (cf. Proposition 4.6).   The following theorem  shows that they are the same if and only if $H^d_\m(M)$ satisfies the property (*).  A characterization of the property (*) for $H^d_\m(M)$ in term of systems of
parameters is also given.

\proclaim{Theorem 3.4}  The following statements are equivalent:
\medskip
(i) $H^d_\m(M)$ satisfies the property (*).
\smallskip
(ii) $\Usupp M= \{\widehat\p \cap R :\ \widehat \p\in\Usupp_{\R}
\widehat M\}.$
\smallskip
(iii) For every sequence $x_1,\ldots ,x_d$ of elements in $\m$,
$(x_1,\ldots ,x_d)$ is a system of parameters of $H^d_\m(M)$ if
and only if it is a system of parameters of $M/U_M(0).$
\endproclaim

\demo{Proof} (i) $\Leftrightarrow$ (ii). We get by Lemma 3.2 that
$V(\Ann H^d_\m(M))=\Usupp M$ and
$V(\Ann_{\R}H^d_\m(M))=\Usupp_{\R}\widehat M.$ Therefore the condition
(ii) is equivalent to the condition
$$V(\Ann H^d_\m(M))=\{\widehat\p \cap R\ :\ \widehat\p\in
V(\Ann_{\R}H^d_\m(M))\}.$$ 
So,  our claim  follows by Proposition 2.2.

\noindent (i) $\Rightarrow$ (iii). Let $(x_1,\ldots, x_d)$ be a
system of parameters of $H^d_\m(M).$ Let $I$ be the ideal
generated by $x_1,\ldots ,x_d.$
 For each prime ideal $\p$  of $R$  containing $I+\Ann H^d_\m (M),$
 we have by (i) that
$$\p = \Ann (0:_{H^d_\m(M)}\p )\supseteq \Ann (0:_{H^d_\m(M)}I).$$
 Therefore
$$ \rad \big(I+\Ann H^d_\m (M)\big)=\bigcap_{\p\supseteq I+
\Ann H^d_\m (M)}\p\supseteq \rad \big(\Ann
(0:_{H^d_\m(M)}I)\big).$$
 Hence $\rad \big(I+\Ann H^d_\m
(M)\big)=\rad \big(\Ann (0:_{H^d_\m(M)}I)\big).$ Since $(x_1,\ldots ,x_d)$ is a
system of parameters of $H^d_\m(M),$ the length of
($0:_{H^d_\m(M)}I$) is  finite.  So, we get by the last equality that $I+\Ann H^d_\m (M)$ is an $\m-$primary ideal. Since $\rad\big(\Ann H^d_\m (M)\big)=\rad\big(\Ann (M/U_M(0))\big),$ (cf. Lemma 3.2), the ideal $I+\Ann (M/U_M(0))$ is
 $\m-$primary. Hence $(x_1,\ldots ,x_d)$ is a system
of parameters of $M/U_M(0).$ Conversely, assume that $(x_1,\ldots
,x_d)$ is a system of parameters of $M/U_M(0).$ Then $I+\Ann
(M/U_M(0))$ is $\m-$primary, and hence so is $I+\Ann H^d_\m (M).$
Therefore $\ell (0:_{H^d_\m(M)}I)<\infty$, i.e. $(x_1,\ldots
,x_d)$ is a system of parameters of $H^d_\m(M).$
\smallskip

\noindent (iii)$\Rightarrow$ (i). Let $\p\in V(\Ann H^d_\m(M)).$
Assume that $\N-dim (0:_{H^d_\m(M)}\p)=d-r.$  By [TZ, Proposition 2.10], there exist  $x_1,\ldots ,x_r\in \p$ which form a  part of a
system of parameters of $H^d_\m(M)$ in $\p ,$ and it is clear that this part of a system of parameters is maximal. Let
$$0:_{H^d_\m(M)}(x_1,\ldots ,x_r)R=A_1+\ldots +A_n$$
  be a minimal secondary representation of $0:_{H^d_\m(M)}(x_1,\ldots ,x_r)R,$
   where $ A_i$ is $\q_i-$secondary.  For each element $y\in \m$, note that $y$  is a parameter element of $0:_{H^d_\m(M)}(x_1,\ldots ,x_r)R$ if and only if $y\not\in\q_i$ for all $i$ satisfying $\N-dim A_i=d-r$ (cf. [TZ, Lemma 2.14]).  Since $(x_1,\ldots ,x_r)$ is a maximal  part of a system of parameters of $H^d_\m(M)$ in $\p ,$ we have 
$$\p\subseteq\bigcup_{\N-dim A_i=d-r}\q_i$$
 and hence $\p\subseteq \q_i$ for some $i$ satisfying $\N-dim A_i =d-r.$
  From the hypothesis (iii), we can check  that  $(x_1,\ldots ,x_r)$
  is a maximal part of a system of parameters of $M/U_M(0)$ in $\p .$ So,  there exists a  prime
ideal $\q\in \Ass \big(M/U_M(0)/(x_1,\ldots ,x_r)M/U_M(0)\big)$ such that $\dim R/\q =d-r$ and $\p\subseteq\q$.  Since  $\p\in\Supp
\big(M/U_M(0)/(x_1,\ldots ,x_r)M/U_M(0)\big),$ it follows that $\p=\q .$  Hence $\dim R/\p =d-r.$ Since $A_i$ is $\q_i-$secondary, we have by Lemma 2.4,(i) that
 $\N-dim A_i\leqslant \dim R/\q_i.$  Because $\p\subseteq \q_i,$ we have
$$d-r=\N-dim A_i\leqslant \dim R/\q_i\leqslant \dim R/\p =d-r.$$
Hence $\p =\q_i$ and hence $\p\in\Att (0:_{H^d_\m(M)}(x_1,\ldots ,x_r)R).$
 Therefore there exists a prime ideal  $\widehat\p\in\Att_{\widehat R}
  (0:_{H^d_\m(M)}(x_1,\ldots ,x_r)R)$ such that $\widehat\p \cap R=\p.$ This implies that
$$\p\subseteq \Ann  (0:_{H^d_\m(M)}\p )\subseteq \Ann_{\widehat R} (0:_{H^d_\m(M)}\widehat\p)\cap
R=\widehat\p \cap R=\p .$$  Thus $\Ann  (0:_{H^d_\m(M)}\p )=\p$ as required. \hfill\ \qed
\enddemo

\bigskip

\noindent{\bf 4. The catenaricity of $\Usupp M$}

 \bigskip

 We say that  $\Supp M$ is {\it catenary} if for any prime ideals $\p ,\q\in\Supp M$ with  $\p\subset\q,$ all saturated chains of prime ideals starting from $\p$ and ending at $\q$ have the same length.  It is clear that $\Supp M$ is catenary if and only if the ring $R/\Ann M$ is catenary.    Therefore, $\Supp M$  is catenary and $\dim R/\p =d$ for all minimal prime ideals $\p\in\Ass M$ if and  only if $\dim R/\p +\dim M_{\p} =d$ for all $\p\in\Supp M,$ cf. Remark 18.5 of  [HIO]. In particular,  since $\dim R/\p=d$ for all $\p\in\Ass M/U_M(0),$ the umixed support $\Usupp M=\Supp M/U_M(0)$ is catenary if and only
 if $\dim R/\p +\dim M_{\p}=d$ for all $\p\in\Usupp M.$
\smallskip

  The following theorem is the main result of this section, which  shows that the property (*)
  for $H^d_\m(M)$ is equivalent to the catenaricity of $\Usupp M$.

\proclaim{Theorem 4.1} The following statements are equivalent:

(i) $\Usupp M$ is catenary.

(ii)  $H^d_\m(M)$  satisfies the property (*).
\endproclaim

Before proving Theorem 4.1, we need the following lemmas.

\proclaim{Lemma 4.2} Assume that $R$ is complete with respect to the 
$\m-$adic topology and $M$ a finitely generated $R-$module such
that $\dim R/\p=d$ for all $\p\in\Ass M.$ Then $\dim R/\p=d-r$ for
any part of system of parameters $(x_1,\ldots ,x_r)$ of $M$ and
any minimal associated prime ideal $\p$ of $M/(x_1,\ldots ,x_r)M.$
\endproclaim
\demo{Proof} 
As $(x_1,\ldots ,x_r)$ is a part of a system of parameters  of $M$, we have 
$$\dim (R/\Ann M+(x_1, \ldots, x_r)R)=\dim (M/(x_1, \ldots, x_r)M)=d-r.$$
Moreover $\p$ is a minimal prime divisor of $\Ann M+(x_1, \ldots, x_r)R$, so that $\dim R/\p$ is at most $d-r.$ There is a minimal prime divisor $\q$ of $\Ann M$ which is contained in $\p$. As $\q$ belongs to $\Ass_R(M)$ it follows from our assumptions that $\dim R/\q=d.$ Moreover $\p$ is a minimal prime divisor of $\q+(x_1, \ldots, x_r)R$, so that $\docao (\p/\q)$ does not exceed $r$ (cf. [Mat,  Theorem 18]).  As $R/\q$ is catenary of dimention $d$, it follows that $\dim R/\p=d-\docao (\p/\q)$ is at least $d-r$. Therefore $\dim R/\p=d-r.$ \hfill\ \qed
\enddemo
\proclaim{Lemma 4.3} Let $\p\in V(\Ann  H^d_\m(M))$ such that
$\dim M_\p +\dim R/\p =d.$ Then $\Ann (0:_{H^d_\m(M)}\p )=\p .$
 \endproclaim

\demo{Proof}  Let $\p\supseteq\Ann  H^d_\m(M)$ be a prime ideal
such that $\dim M_\p +\dim R/\p =d.$  Set $\dim R/\p =d-r.$ It
follows by the hypothesis that $\dim M_\p =r.$ Therefore  there
exists a prime ideal $\q\in\Ass M$ such that $\q\subseteq\p $ and
$\docao (\p /\q )=r.$ Since
$$\dim R/\q \geq \dim R/\p +\docao (\p/\q )=d,$$ 
it follows that $\dim R/\q =d.$ It should be noted that $\dim \R /\p\R =\dim R/\p
=d-r.$ So,  there exists a prime ideal $\widehat\p\in\Ass_{\R}\widehat
R/\p\widehat R$ such that $\dim\widehat R/\widehat\p =d-r.$ Since
$\widehat\p\in\Ass_{\R} \widehat R/\p\widehat R,$ we get
$\widehat\p\cap R\in\Ass R/\p ,$ i.e. $\widehat\p\cap R=\p .$ Note
that the natural map $R\lr \R$ is faithfully flat, and therefore the
going down theorem holds (see [Mat, Theorem 4]). So, there
exists a prime ideal $\widehat\q\in\Spec\widehat R$ such that
$\widehat\q\cap R=\q , \ \widehat\q\subseteq\widehat\p$ and
$\docao (\widehat\p /\widehat\q )\geq r.$ These facts imply that
$$d=\dim R/\q =\dim \widehat R/\q\widehat R
\geq \dim \R/\widehat\q =\dim\widehat R/\widehat\p +\docao
(\widehat\p /\widehat\q )\geq d-r +r=d.$$ Hence $\dim
\R/\widehat\q =d.$ Moreover, since the natural homomorphism
$R_\q\lr\widehat R_{\widehat\q}$ is faithfully flat and $M_\q \neq
0,$ we have
$$M_\q\otimes_{R_\q}\widehat R_{\widehat\q}\cong \widehat M_{\widehat\q}\neq 0.$$
Hence  $\widehat\q\in\Supp_{\R}\widehat M.$ Since $\dim \R/\widehat\q
=d$ and $\widehat\p\supseteq\widehat\q ,$ we have
$\widehat\p\supseteq\Ann_{\R} H^d_\m(M).$ Therefore we get by the
Matlis duality that $\Ann_{\widehat R}(0:_{H^d_\m(M)}\widehat\p
)=\widehat\p .$ Now we have
$$\p\subseteq\Ann  (0:_{H^d_\m(M)}\p )\subseteq \Ann_{\widehat R} (0:_{H^d_\m(M)}\widehat\p )\cap R
=\widehat\p\cap R =\p .$$ Thus $\Ann  (0:_{H^d_\m(M)}\p )=\p .$
\hfill\ \qed
\enddemo

\noindent{\it Proof of Theorem 4.1.}  (i) $\Rightarrow$ (ii)
follows by Lemma 4.3 and Lemma 3.2.

\noindent (ii) $\Rightarrow$ (i).   Let $\p\in\Usupp M.$ We need to show that $\dim R/\p +\dim M_{\p}=d$.  The case $\p=\m$ is trivial. Assume that $\p\neq\m .$ Let $\dim R/\p =d-r.$
Then it is enough to prove  $\dim M_{\p}=r.$ Since
$\p\supseteq\Ann  M/U_M(0),$ we have
$$\dim \Big(M/U_M(0)\big/\p (M/U_M(0))\Big)=\dim R/\p =d-r.$$ So,
there exists a maximal part of a system of parameters $(x_1,\ldots
,x_r)$ of $M/U_M(0)$ in $\p .$ Since $\p\in\Usupp M,$ there exists
by Theorem 3.4, (i)$\Leftrightarrow$(ii)  a prime ideal
$\widehat\p\in\Usupp_{\R} \widehat M$ such that $\widehat\p \cap R=\p
.$ Set $\widehat M_1=\widehat M/U_{\widehat M}(0).$ Since
$(x_1,\ldots ,x_r)$ is a part of a system of parameters of
$M/U_M(0),$ it is a part of a system of parameters of the
$\m-$adic completion $\widehat{M/U_M(0)}$ of $M/U_M(0)$. Because
$\widehat M_1$ is a quotient of $\widehat{M/U_M(0)}$ and $\dim
\widehat M_1 =\dim \widehat{M/U_M(0)},$ it follows that
$(x_1,\ldots ,x_r)$ is a part of a system of parameters of
$\widehat M_1.$  Note that $\widehat\p\in\Supp_{\R} \widehat
M_1/(x_1,\ldots ,x_{r-1})\widehat M_1.$ Therefore
$\widehat\p\supseteq \widehat\p_1$ for some minimal prime ideal
$\widehat\p_1\in\Supp_{\R}\widehat M_1/(x_1,\ldots ,x_{r-1})\widehat
M_1.$ Since $x_r$ is a parameter element of $\widehat
M_1/(x_1,\ldots ,x_{r-1})\widehat M_1,$ we get by Lemma 4.2 that
$x_r\notin\widehat\p_1.$  Set $\p_1=\widehat\p_1\cap R.$ Then
$x_r\notin\p_1$. Therefore $\p\supset \p_1$ and $\p\neq \p_1.$ By
the same arguments, there exists a minimal prime ideal
$\widehat\p_2\in\Supp_{\R}\widehat M_1/(x_1,\ldots ,x_{r-2})\widehat
M_1$ such that $\widehat \p_1\supseteq\widehat \p_2.$ Set
$\p_2=\widehat\p_2\cap \R.$ Then $\p_1\supset\p_2$ and $\p_1\neq
\p_2$ since $x_{r-1}\in\p_1\setminus\p_2.$ Continue the above
process,  after $r$ steps,  we get a chain $\p\supset \p_1\supset
\p_2\ldots \supset \p_r$ of prime ideals containing $\Ann M$ such
that $\p_i\neq \p_{i+1}$ for all $i=1,\ldots r-1.$ Therefore $\dim
M_\p =r.$ \hfill\ \qed
\medskip

\noindent{\bf Proof of the main theorem.} The equivalences between (i),
(iii) and (iv) are proved by Theorem 3.4. The equivalence between
(i) and (ii) is proved by Theorem 4.1. \hfill\ \qed
\medskip

\noindent{\bf Remark 4.4.} The catenaricity of $\Usupp M$ is
equivalent to the property (*) for the top local cohomology module
$H^d_\m(M)$ (see Theorem 4.1), but it is not related to the
property (*) of other local cohomology modules of $M.$ In fact,
let $R$ be the Noetherian local domain constructed by Ferrand  and
M. Raynaund [FR] of dimension $2$ such that the $\m-$adic
completion has an associated prime of dimension $1.$ It is clear
that $\Supp R=\Usupp R$ is catenary, but $H^1_\m(R)$ does
not satisfy the property (*).
\medskip

 Let $0=M_0\subset M_1\subset \ldots \subset M_t=M$ be  the filtration of submodules of $M$, where  $M_{i-1}$ is the largest submodule of $M_i$ of dimension less than $\dim M_i$ for $i=1,\ldots ,t.$ Such a filtration always exists, and it is unique. We call this filtration to be {\it  the dimension filtration} of $M$ (cf. [CN2]).  Let $\dim M_i=d_i$ for $i=1,\ldots ,t.$ Then it is easy to check that $$\Supp M=\underset{i=1,\ldots ,t}\to{\bigcup}\Supp M_i/M_{i-1}.$$
For each $i=1,\ldots ,t,$ it should be noted that $\dim R/\p =d_i$ for all $\p\in\Ass M_i/M_{i-1}.$ Therefore  we obtain by Theorem 4.1 the following result.

\proclaim{Corollary 4.5} $\Supp M$ is catenary if and only if $H^{d_i}_\m (M_i/M_{i-1})$ satisfies the property (*) for all $i=1,\ldots ,t.$
\endproclaim

Now we examine some non-catenary Noetherian local domains. Note
that any domain of dimension $2$ is catenary, but there exist non-catenary Noetherian local domains of dimension $d$ for any $d\geq 3$  (cf. [B, (8)]).

 \proclaim{Proposition 4.6} Let $R$ be a non-catenary Noetherian local domain of dimension $3.$  Set
$$\align &U=\{\p\in\Spec R : \dim R/\p +\docao \p =2\};\\
&V=\{\p\in\Spec R : \dim R/\p +\docao \p =3\}.\endalign$$
  Then the following statements are true
\medskip

(i) $\Usupp R=\Spec R=U\cup V$ and  $U, V\neq \emptyset .$

\smallskip

(ii) $\Ann (0:_{H^3_\m(R)}\p )=\p$ for all $\p\in V.$ But $\Ann
(0:_{H^3_\m(R)}\p )\neq \p$ for all $\p\in U$.
\smallskip
(iii) For each $\p\in V,$ there always exists $\widehat\p
\in\Supp\R/U_{\R}(0)$ such that $\widehat\p\cap R=\p .$ But for
each $\p\in U,$ there does not exist $\widehat\p \in\Supp
\R/U_{\R}(0)$ such that $\widehat\p\cap R=\p .$
\smallskip
(iv) $\N-dim H^2_\m(R) =2$ and $\dim R/\Ann  H^2_\m (R)=3.$

\endproclaim

\demo{Proof}  (i). This is clear  since $R$ is a non-catenary domain.

\noindent (ii). It follows by the proof of Theorem 4.1 that $\Ann
(0:_{H^3_\m(R)}\p )\neq \p$ for all $\p\in U$ and $\Ann
(0:_{H^3_\m(R)}\p )=\p$ for all $\p\in V.$
\smallskip
\noindent (iii). This  follows by (ii) and by the proof of Theorem 3.4.

\smallskip

\noindent (iv). Let $\p\in U.$ Then $\dim R/\p =1.$
 Let $\widehat\p\in\Spec\widehat R$ such that $\widehat\p\cap R=\p .$
 Then $\dim\widehat R/\widehat\p =1.$  It follows by (iii) that
 $\widehat\p\not\supseteq \Ann_{\R} H^3_\m(R).$
  Moreover,  $\docao \widehat\p \geq \docao\p =1$
   by the going down theorem [Mat, Theorem 4].
    Therefore  there exists $\widehat\q\in\Ass\widehat R$
    such that $\widehat\q\subset\widehat\p$ and  $\widehat\q\neq \widehat\p .$
     Hence  $\dim\widehat R/\widehat\q \geq 2.$ Since $\widehat\p\not\supseteq \Ann_{\R} H^3_\m(R),$
     it follows that  $\dim\widehat R/\widehat\q
     =2.$  So we have by [BS, Corollary 11.3.3]
     that  $\widehat\q\in\Att_{\widehat R} H^2_\m(R)$
     and hence $\widehat\q\supseteq \Ann_{\R}H^2_\m(R).$
     Therefore $\N-dim H^2_\m(R)=\dim \R/\Ann_{\R}H^2_\m(R)\geq 2.$ Note that $\N-dim H^2_\m(R)\leqslant 2$
      by [CN1, Theorem 3.1]. So $\N-dim H^2_\m(R)= 2.$
   Since $\widehat\q\in\Att_{\widehat R}  H^2_\m(R)
      \cap \Ass \widehat R,$ we have $$\widehat\q \cap R\in\Att
      H^2_\m(R)\cap \Ass R.$$ Since $R$ is a domain, we have $\widehat\q \cap R=0.$
      It follows that $0=\Ann H^2_\m(R).$ Thus $\dim R/\Ann H^2_\m(R)
      =3.$\ \hfill\ \qed
      \enddemo

\medskip
\noindent{\bf Acknowledgment.}  We wish to express our gratitude to the refree  for his/her  useful suggestions and, especially, the shorter proof of  Lemma 4.2.  

\bigskip
 \noindent{\bf References}
\medskip
\item {[B]} M. Brodmann, {\it A particular class of regular domains},
 J. Algebra, {\bf 54}, (1978), 366-373.

\item {[BS]} M.  Brodmann   and R. Y. Sharp, ``Local cohomology:
an algebraic introduction with geometric applications",
Cambridge University Press, 1998.

\item {[CN1]}  N. T. Cuong and L. T. Nhan, {\it On the Noetherian
dimension of
 Artinian modules}, Vietnam J. Maths.,  {\bf(2)30} (2002), 121-130.

\item{[CN2]} N. T. Cuong and L. T. Nhan,  On pseudo Cohen-Macaulay
and pseudo generalized Cohen-Macaulay modules, {\it J. Algebra,}
{\bf 267} (2003), 156-177.

 \item{[FR]} D. Ferrand  and M. Raynaund,  {\it Fibres formelles d'un
anneau local Noetherian,} Ann. Sci. E'cole Norm. Sup.,}  {\bf
(4)3} (1970), 295-311.

\item{[HIO]} M. Herrmann, S. Ikeda and U. Orbanz,
``Equimultiplicity and Blowing up", Springer - Verlag, 1988.

\item {[K1]} D. Kirby, {\it Artinian modules and Hilbert
polynomials,} Quart. J. Math. Oxford, {\bf (2)24} (1973), 47-57.

 \item {[K2]}  D. Kirby, {\it  Dimension and length of Artinian modules},
  Quart. J. Math. Oxford., {\bf (2)41} (1990), 419-429.

  \item{[Mac]} I. G. Macdonald, {\it  Secondary representation of modules
   over a commutative ring,}  Symposia Mathematica, {\bf 11} (1973), 23-43.

\item {[Mat]}  H. Matsumura, ``Commutative Algebra", Second Edition (Benjamin, 1980).

\item{[Na]} M.  Nagata,  {\it Local rings,} Interscience, New York, 1962.

 \item {[R]} R. N. Roberts, {\it Krull dimension for Artinian modules over
quasi local commutative rings,}
 Quart. J. Math. Oxford, {\bf (2)26} (1975), 269-273.

\item{[Sh]} R. Y. Sharp, {\it A method for the study of
Artinian modules with an application to asymptotic behaviour,}  Commutative
Algebra  (Math. Sciences Research Inst. Publ. No. 15, Spinger Verlag), (1989), 443-465.

\item{[TZ]} Z. Tang  and H. Zakeri,  {\it Co-Cohen-Macaulay modules and
modules of generalized fractions},  Comm. Algebra.,  {\bf (6)22}
(1994), 2173-2204.

\end